\documentclass[aps, pre, twocolumn, groupedaddress, floatfix]{revtex4-1}

\RequirePackage{color}
\RequirePackage[colorlinks, urlcolor=my-blue,linkcolor=my-red,citecolor=my-green]{hyperref}
\definecolor{my-blue}{rgb}{0.0,0.0,0.6}
\definecolor{my-red}{rgb}{0.5,0.0,0.0}
\definecolor{my-green}{rgb}{0.0,0.5,0.0}
\definecolor{nicos-red}{rgb}{0.75,0.0,0.0}
\definecolor{light-gray}{gray}{0.6}
\definecolor{really-light-gray}{gray}{0.8}
\definecolor{sussexg}{rgb}{0.0,0.5,0.5}
\definecolor{sussexp}{rgb}{0.5,0.0,0.5}

\usepackage{nicefrac, natbib, amsmath, amssymb, graphicx,tikz}
\usepackage{graphics}
\newcommand{\be}{\begin{equation}}
\newcommand{\ee}{\end{equation}}

\def\cS{\mathcal{S}}

\def\cL{\mathcal{L}}

\def\bE{\mathbb{E}}
\def\bN{\mathbb{N}}
\def\bP{\mathbb{P}}

   \def\N{\bN}

\def\e{\varepsilon}

\usepackage[british]{babel}

\def\m1{\mathbf{1}}

\def\E{\bE}
\def\P{\bP} %% environment measure 
\def\tmix{t_{\text{mix}}}

\definecolor{darkgreen}{rgb}{0.0,0.5,0.0}
\definecolor{darkblue}{rgb}{0.0,0.0,0.3}
\definecolor{nicosred}{rgb}{0.65,0.1,0.1}
\definecolor{light-gray}{gray}{0.7}

\begin{document}


\title{\textcolor{black}{Solvable non-Markovian dynamic network}}

\author{Nicos Georgiou}
\email[]{N.Georgiou@sussex.ac.uk}
%\thanks{}
\affiliation{School of Mathematics and Physical Sciences, University of Sussex}

\author{Istvan Z. Kiss}
\email[]{I.Z.Kiss@sussex.ac.uk}
%\thanks{}
\affiliation{School of Mathematics and Physical Sciences, University of Sussex}

\author{Enrico Scalas}
\email[]{E.Scalas@sussex.ac.uk}
%\thanks{}
\affiliation{School of Mathematics and Physical Sciences, University of Sussex}


\date{\today}

\begin{abstract}
Non-Markovian processes are widespread in natural and human-made systems, yet 
explicit modelling and analysis of such systems is underdeveloped. We consider a non-Markovian
dynamic network with random link activation and deletion (RLAD) and heavy tailed Mittag-Leffler distribution for the inter-event times.
%thus considerably slowing down the corresponding Markovian dynamics and study the system far from equilibrium.  
%\textcolor{black}{This
%distribution interpolates between a stretched exponential for small inter-event times and a power law of index $\beta <1$ for
%large inter-event times.}
We derive an analytically and computationally tractable system of Kolmogorov-like forward equations utilising the Caputo derivative for the probability of having a given number of active links in the network and solve them.  
\textcolor{black}{Simulations for the RLAD are also studied for power-law inter-event times and we show
excellent agreement with the Mittag-Leffler
model. This agreement holds even when the RLAD network dynamics is coupled with the susceptible-infected-susceptible (SIS) spreading
dynamics. Thus, the analytically solvable Mittag-Leffler model provides an excellent approximation 
to the case when the network dynamics is characterised by power-law distributed inter-event times. We further discuss possible generalizations of our result.}

\end{abstract}



\maketitle

% References should be done using the \cite, \ref, and \label commands%
%\section{Introduction}

\section{Introduction}

Non-Poisson temporal statistics where time intervals between isolated, consecutive actions
are typically not exponentially distributed, seem to be the norm rather than the exception for many systems:
For example, period of infectiousness \cite{lloyd2001realistic}, inter-order and inter-trade durations in financial markets \cite{scalas2006durations},
socio-networks, including emails \cite{eckmann2004entropy,malmgren2008poissonian}, phone calls \cite{jiang2013calling}, or individual-to-individuals contacts being fluid 
\cite{schneider2013unravelling,moinet2014burstiness}.
The absence of the robust
tools and mathematical machinery of Markovian theory is the source of
many challenges in modelling and analysis of non-Markovian systems.
The burst in research activity that successfully combines
networks and non-Markovian processes stems from
the need to develop more realistic models
and new analytical tools. Notable examples include studying non-Poisson dynamics  of networks 
\cite{hoffmann2012generalized} and non-Markovian epidemics on networks 
\cite{min2011spreading,van2013non,jo2014analytically}. 

The non-Markovian property is particularly pervasive when considering the dynamics of time-evolving 
networks, be it with fast or slow timescale \cite{holme2012temporal,perra2012activity}. 
%Given that dynamic networks already pose
%significant modelling and analysis challenges due to their richer state space and behaviour, 
%combining it with, or studying it in, the non-Markovian framework leads to a significant increase
%in complexity and potential lack of analytical tractability. 
%Therefore, deriving simple, solvable paradigm/mechanistic models
%in this context, allows for important progress in understanding the nature of the problem. Moreover, it develops new
%mathematical tools and methods for analysis and increases our understanding of the true implications for complex systems.
%
Deriving simple, solvable
paradigm models can facilitate progress 
%in understanding the nature and the challenges of
%the problem. More importantly, it leads {\color{black} to} the 
in developing new mathematical tools and methods for analysis and increases our understanding
of the true implications of non-Markovianity for complex systems.
{Empirically, it turns out that many inter-event distributions have power-law tails
(see \cite{vajna2013} and references therein). Therefore, it is also necessary to develop methods 
able to deal with such distributions.}

It is now widely accepted that human contact patterns are highly dynamic and may evolve concurrently with
an epidemic; many Markovian models for this setup exists \cite{gross2009adaptive,marceau2010adaptive,kiss2012modelling}.  
Here, we take the next step and consider a dynamic network with non-exponential waiting times with consecutive 
updates which are either link activation and deletion \cite{kiss2012modelling}. 
As a first step in the rigorous analysis
of networks with Non-Markovian dynamics, we consider a random link activation-deletion (RLAD) 
model that naturally leads to a stochastically evolving network \cite{raberto2011graphs, kiss2012modelling}. This model amounts to
considering undirected and unweighted networks, where an event consists of selecting a link at random, 
independently of whether present or not, followed by its activation, if the link is absent, or deletion if the link is active. 
Such operations are separated by inter-event times sampled from the Mittag-Leffler distribution, that allows for analytical tractability.
This exactly solvable model of non-Markovian network dynamics is an important special case 
 of a more general theory for non-Markovian 
processes outlined in \cite{haenggi1981old}, and 
it is related to recent outstanding developments in probability theory \cite{orsingherpolito,meerschaerttoaldo}. 
Indeed, we provide a bottom-up derivation for the master equation of some 
fractional birth and death processes 
in a finite capacity system, 
introduced in \cite{orsingherpolito}.
This allows us to compute theoretically the exact distribution of the total number of links in the network at any time
and its large-time limit. 
%and the consequences for a simple $SIS$ epidemic on this network. 
%we explain why the choice for a distribution should not be considered restrictive;
We demonstrate the power of the 
analytical model by comparing it with simulations using more widely-used power-law distributed times.
The rigorous analysis of this model, including explicit expressions for the distribution of the number of links in the network for $t\ge 0$, 
is followed by considering a Markovian $SIS$ epidemic on our non-Markovian dynamic network. \textcolor{black}{Finally, we briefly discuss the generalization of our method to general Markov chains with random state changes occurring according to a generic renewal process}.
	
%the dynamic network will be coupled with a  model to illustrate the impact of non-Marovian network dynamics on epidemic dyanics.
%For example, in order to illustrate how this theory can be used with applications to epidemiology, 
%we simulate an SIS epidemic on the RLAD network. At that point we assume the initial condition of the network
%to be the complete graph since that  is the one that makes sense from an epidemiology perspective. 
%As such, we begin with 
%simplified models that are analytically tractable in order to start building a theory 
%that encompasses more general models and eventually 
%deals with the modeling and statistical problems discussed earlier.

\section{An exactly solvable model}

\subsection{Basic ingredients}

Consider an arbitrary graph on $N$ nodes as 
an initial state of the dynamics. 
We are interested in the number of (unique undirected) links in the network at a given time $t$. 
We denote this number by $X(t)$, and it takes values in $\cS=\{ 0, 1, \ldots, M\}$ 
where $M = N(N-1)/2$, the maximal possible number of links. 
The time periods where $X(t)$ remains constant are called sojourn times or inter-event times. 
We assume that sojourn times $\{ T_i\}_{i \ge 1}$ 
are drawn independently from the family of Mittag-Leffler distributions with parameter (or order)
$\beta \in (0,1)$  \cite{frac-calc-henry}. Their cumulative distribution function (c.d.f.) \ is indexed by this $\beta$ and it is given by 
	\be\label{eq:m-l:cdf}
		F^{(\beta)}_{T}(t)= \bP\{ T \le t \} = 1 - E_{\beta}(-t^{\beta}).
	\ee
Here $E_{\beta}(z)$ is the Mittag-Leffler function, defined by 
	\be \label{eq:mlf}
		E_{\beta}(z) = \sum_{n=0}^{\infty} \frac{z^n}{\Gamma(1 + \beta n)}.
	\ee
$E_{\beta}$ is entire for all $\beta > 0$. At $\beta = 0$ the series converges uniformly only on a disc 
of radius 1, though the function can be extended analytically on $\mathbb{C} \smallsetminus \{1\}$. 
Equations \eqref{eq:m-l:cdf}, \eqref{eq:mlf} define a proper c.d.f. only when $\beta \in (0,1]$. \textcolor{black}{This is equivalent to the claim that, for $\beta \in (0,1]$, $E_\beta (-t^\beta)$ is completely monotone. A $C^\infty[0,\infty)$ function $f(t)$ is completely monotone if $(-1)^n d^n f(t)/dt^n \geq 0$ for all non-negative integer $n$ and all $t>0$. Now, Mainardi and Gorenflo \cite{mainardigorenflo} proved that, for $\beta \in (0,1)$, $E_\beta(-t^\beta)$ can be written as a mixture of exponential distributions given that
\begin{equation}
E_\beta (-t^\beta) = \int_0^\infty \exp(-rt) K_\beta (r) \, dr,
\end{equation}
where
\begin{equation}
K_\beta (r) = \frac{1}{\pi} \frac{r^{\beta-1} \sin(\beta \pi)}{r^{2\beta} + 2 r^\beta \cos(\beta \pi) + 1},
\end{equation}
and
\begin{equation}
\int_0^\infty K_\beta (r) \, dr =1.
\end{equation}
Therefore, complete monotonicity of $E_\beta (-t^\beta)$ is an immediate corollary of Bernstein's theorem \cite{bernstein,schillingsongvondracek}. A direct proof that $E_\beta (-x)$ is completely monotone can be found in reference \cite{pollard}.}
When $0<\beta <1$ these distributions are heavy-tailed with infinite mean
while at $\beta = 1$, $T$ is  mean $1$, exponentially distributed.
This family of distributions interpolates between a stretched exponential for small $t$ and a power-law
for large $t$ \cite{mainardigorenflo}. Namely, one has
\begin{eqnarray}
E_\beta (-t^\beta) & \simeq & \exp(-t^\beta/\Gamma(1+\beta)), \,\, t  \ll 1, \nonumber \\
E_\beta (-t^\beta) & \sim & \frac{\sin(\beta \pi)}{\pi} \frac{\Gamma(\beta)}{t^\beta}, \, \, t \to \infty \label{eq:ref:stupid}.
\end{eqnarray}
Therefore, the use of these distributions is more general than it might seem at a first glance.  
A word of notational caution: Here $\beta$ is the order of the polynomial decay of the survival function, but most commonly power-law distributions are identified by the order of decay of their densities, which in our case is $1 + \beta \in (1,2)$. 

Mittag-Leffler sojourn times lead to 
a simpler analytical treatment of non-Markovianity in the presence of extreme power-law tails than its cognate Pareto distribution. 
However, we do explain below how the theoretical framework developed here can be used to approximate the behaviour 
of non exactly-solvable systems with Pareto power-law distribution, as it is most commonly used. 
For this we must introduce a scaling parameter (time change)
$\gamma > 0$ for the waiting times: We say that a random variable $T$ is Mittag-Leffler$_\gamma(\beta)$ distributed if and only if 
\be\label{eq:MLgamma}
F^{(\beta,\gamma)}_{T}(t)= \bP\{ T \le t \} = 1 - E_{\beta}(- (t/\gamma)^{\beta}).
\ee  
For $\gamma=1$ the c.d.f.\! is reduced to that of equation \eqref{eq:m-l:cdf} and we see that 
$T$ is Mittag-Leffler$_1(\beta)$ if and only if $\gamma T$ is Mittag-Leffler$_\gamma(\beta)$. 

For the rigorous derivation of the evolution equations, we restrict for clarity to the $\gamma=1$ case and remark how 
the equations behave with the extra scaling later.
Fix a parameter $\beta \in (0,1) $. The network evolves in a semi-Markov way: 
Let $T_1, T_2, \ldots$ be  independent Mittag-Leffler$(\beta)$ times and define 
the partial sum \textcolor{black}{
\begin{equation}
S_n = \sum_{k=1}^n T_i, \; \; n\ge 1. 
\end{equation}
}The sequence $S_1, S_2, \ldots $ denotes the event times at 
which the state of the network $X(t)$ attempts to change. 
A change in the state means
an undirected link is either deleted or activated. For extra flexibility, the model is introduced with 
an extra delay parameter $\alpha \in [0,1)$, so that if $\alpha \neq 0$  allows the active links to remain 
unchanged even if there is an attempt of a change.   

It is useful to define the embedded Markov chain for the number of links in the network, 
$X_{n}, n\ge 1$, with state space 
$\cS $. 
Initially $X_0 =i$, as we start with $i$ present links and the number of links in the network increases, remains or 
decreases according to the following transition probabilities
\begin{multline}
\label{eq:trans1}
	q_{k, k-1} =  P_0\{ X_{j+1} = k-1 | X_j = k \} \\
		=  \begin{cases} 
		0, & k =0,\\
		1-\alpha, & k = M,\\
		(1-\alpha)\frac{k}{M}, & \text{ otherwise, }
			\end{cases}
\end{multline}
	\be\label{eq:trans3}
	q_{k,k}=\alpha,
	\ee
	and 
\begin{multline}
\label{eq:trans2}
	q_{k, k+1} = P_0\{ X_{j+1} = k+1 | X_j = k \} \\
		=  \begin{cases} 
		1-\alpha, &k =0,\\
		0,  &k = M,\\
		1- \alpha-\frac{k(1-\alpha)}{M}, &\text{ otherwise. }
			\end{cases}
\end{multline}
In words, at the time of the $i$-th event, we pick a link %$e$ 
uniformly at random 
out of all available links. With probability $\alpha$ nothing changes, otherwise on the event that a change will happen 
in the system,  we delete or add a link in the following way: If the link was active (present) in the network, it is now deleted, 
otherwise it is now activated. Notice that the embedded dynamics are equivalent to the 
$\alpha$-delayed version of the Ehrenfest chain. 

To connect the embedded chain $X_n$ with process $X(t)$, define the counting process 
	\be \label{eq:m-l:count} 
		N_{\beta}(t) = \max\{ n \in \N : S_n \le t \}
	\ee
that gives the number of events up to a finite time horizon $t$. This process is also called a fractional Poisson process. 
Then we have
	\be\label{eq:subordinator}
		X(t) = X_{N_{\beta}(t)} = X_n 1\!\!1\{ S_n \le t < S_{n+1} \},
	\ee
i.e. the state of the process at time $t$ is the same as that 
of the embedded chain after the last event before time $t$ occurred. 

\subsection{Semi-Markov Master Equation} 

All information 
about $X(t)$ is encoded in the pairs 
$\{(X_n, T_n)\}_{n \ge 1}$ which are a discrete-time Markov renewal process, satisfying 
	\begin{align}\label{eq:mrp}
		\P\{ X_{n+1} = j, & T_{n+1} \le t | (X_0, S_0), \ldots, (X_n =i, S_n)\} \notag\\
				&= \P\{ X_{n+1} = j, T_{n+1}\le t | X_n = i \}.
	\end{align}
$X(\cdot)$ is then a semi-Markov process subordinated to $N_\beta (t)$ \cite{raberto2011graphs} and satisfies the forward equations
	\begin{align} 
		p_{i,j}(t) \!= \overline F^{(\beta)}_T(t)\delta_{ij} 
			\!\!+\! \sum_{\ell \in \cS} q_{\ell, j} \!\! \int_0^t \!\! p_{i,\ell}(u)f^{(\beta)}_{T}(t-u)\,du.
			\label{eq:master}
	\end{align}
\textcolor{black}{Incidentally, a semi-Markov process is Markovian if and only if the distribution of $\{T_n\}_{n\geq 1}$ is exponential \cite{cinlar}}.
Above we introduced $p_{i,j}(t) = \P\{X(t) = j | X(0) = i\}$, the tail (complementary cumulative distribution function) 
$\overline F^{(\beta)}_T(t) = 1-F^{(\beta)}_T(t)$ and 
$f^{(\beta)}_{T} (t)$ the Mittag-Leffler density or order $\beta$.
These equations are proved by conditioning on the time of the last event before time $t$. 
	By taking Laplace transforms in \eqref{eq:master}, and using the known Laplace transform 
	of the Mittag-Leffler survival function and probability density function
	\be\label{eq:ml-lap}
		\cL\left(  \overline F^{(\beta)}_T(t) ; s \right) = \frac{s^{\beta-1}}{1 + s^{\beta}}  \text{ and } 
		\cL\left( f^{(\beta)}_T(t) ; s \right) = \frac{1}{1 + s^{\beta}},
	\ee
        followed by some straightforward algebra, the evolution equations for $p_{i,j}(t)$ become (\textcolor{black}{see Appendix A}),
	\begin{align}\label{eq:final0}
		&\frac{d^{\beta} p_{i,j}(t)}{ d\, t ^{\beta}}  = -(1-\alpha)\ p_{i,j}(t) \\
			&\phantom{xxx}+ (1-\alpha)\left(\frac{M -j +1}{M} p_{i, j-1}(t) + \frac{j+1}{M}p_{i,j+1}(t)\right).\notag
	\end{align}
	Similarly, the equations of the boundary terms are
	\be\label{eq:final1}
		\frac{d^{\beta} p_{i,0}(t)}{ d\, t ^{\beta}}  = (1-\alpha)\left(- p_{i,0}(t)+ \frac{1}{M}p_{i,1}(t)\right) \quad %\text{and}
	\ee
	%and
	\be\label{eq:final2}
		\frac{d^{\beta} p_{i,M}(t)}{ d\, t ^{\beta}}  = (1-\alpha)\left(- p_{i,M}(t)+ \frac{1}{M}p_{i,M-1}(t)\right).
	\ee
	Symbol $d^\beta / dt^\beta$ in \eqref{eq:final0}, \eqref{eq:final1}, \eqref{eq:final2}, 
	denotes the $\beta$ {\em fractional Caputo derivative}
\cite{frac-calc-henry} of a function $f(t)$ given by
\[
		\frac{d^{\beta} f(t)}{ d\, t ^{\beta}} = \frac{1}{\Gamma(1-\beta)} \int_{0}^t (t-t')^{-\beta} \frac{d\,f(t')}{dt'}\,dt'. 
\]
	When $\beta = 1$, equations \eqref{eq:final0}, \eqref{eq:final1}, \eqref{eq:final2}
	reduce (as expected) to the standard Kolmogorov equations for the Markovian RLAD \cite{kiss2012modelling}. 
	These equations also explain analytically why $\alpha$ is called the delay parameter.  
	When considering the scaled Mittag-Leffler$_\gamma(\beta)$ times, equation \eqref{eq:final0} becomes 
	\begin{align}\label{eq:final3}
		&\frac{d^{\beta} p^{(\gamma)}_{i,j}(t)}{ d\, t ^{\beta}}  = -\gamma^{-\beta}(1-\alpha)\ p^{(\gamma)}_{i,j}(t) \\
			&\phantom{xx}+ \gamma^{-\beta}(1-\alpha)\left(\frac{M -j +1}{M} p^{(\gamma)}_{i, j-1}(t) +\frac{j+1}{M}p^{(\gamma)}_{i,j+1}(t)\right)\notag
	\end{align}
	and similarly for the boundary equations. Specifically we see, as in the Markovian case, that a scaled 
	sojourn time distribution results in a (fractional) scalar multiple of the forward equations.
	
\subsection{Exact solution}
 
	Equation \eqref{eq:master} gives an analytical way to obtain the fractional equation for the evolution 
	of the transition probabilities, but it is not very useful for computational purposes. 
	Instead, it is fruitful to find the solution of the system of equations \eqref{eq:final0}, \eqref{eq:final1},
	 \eqref{eq:final2} by a simple conditioning argument on the values of $N_{\beta}(t)$ (\textcolor{black}{see Appendix B}) 
	 \be \label{eq:anal:istvan} 
	p_{i,j}(t) = \overline F^{(\beta)}_T(t)\delta_{ij} + \sum_{n=1}^\infty q^{(n)}_{i,j} \P\{N_\beta(t) = n\},
	\ee
	where $q^{(n)}_{i,j}$ are the $n$-step transitions of the embedded discrete Markov chain, namely  the entries  
	of  the $n$-th power  
	of the transition matrix $Q$ defined by equations \eqref{eq:trans1}, \eqref{eq:trans3} and \eqref{eq:trans2}. 
	The distribution of the fractional Poisson process has a simple expression 
	generalising the Poisson distribution \cite{scalas2004}, namely
	\begin{equation}
        \label{eq:m-l:number}
		\mathbb{P}\{N_\beta (t) = n\}= \frac{t^{\beta n}}{n!} E_\beta^{(n)} (-t^\beta),
	\end{equation}
	where $E_\beta^{(n)} (-t^\beta)$ denotes the $n$-th derivative of $E_\beta (z)$ computed for $z = -t^\beta$.
	%If one prefers to avoid the probabilistic argument, 
	Equation \eqref{eq:anal:istvan} can also be verified to satisfy \eqref{eq:final0} using Laplace transforms.
	(\textcolor{black}{see Appendix B}).

\section{Results and applications}

\textcolor{black}{All simulations are event driven, both for dynamic networks
and when this is coupled with epidemic dynamics. Waiting
times for all the possible events are generated from
appropriate distributions. Hence, the next change or an
update is always determined by the smallest waiting
time and the event corresponding to it is executed. This
is then followed by the necessary update of the waiting
times of the events affected by the most recent change.
In reference \cite{boguna}, readers can find an alternative efficient
simulation method which effectively extends the ideas of
the Gillespie algorithm from the Markovian to the non-
Markovian case.}

\subsection{Explicit calculation of $p_{i,j} (t)$}

\begin{figure}[h]
	\begin{center}
	\includegraphics[height=5.5cm,width=8cm]{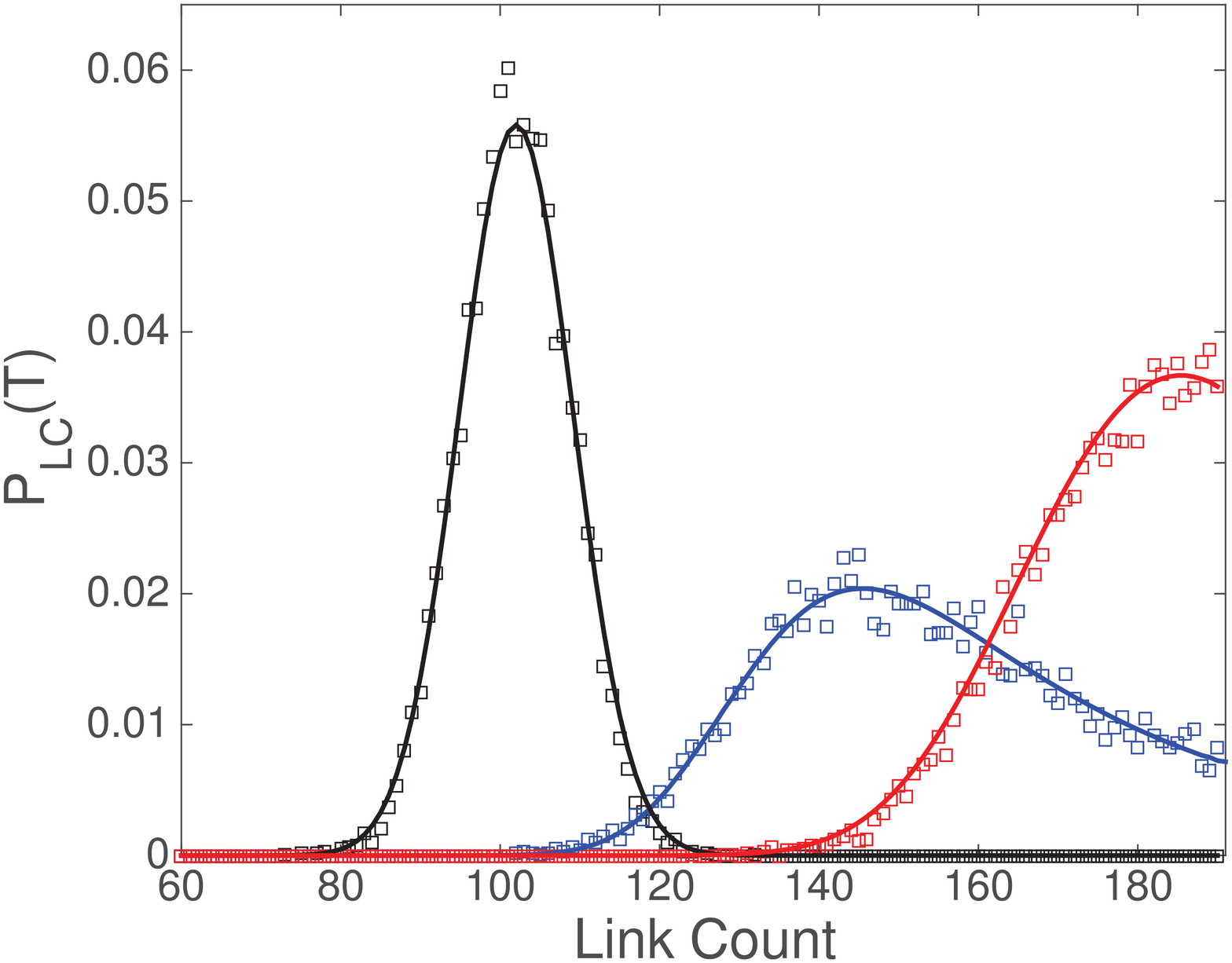}
	\end{center}
	\caption{(Color online) \textcolor{black}{Comparison between Monte Carlo simulations and theory.} The discrete markers are the estimated probabilities $p_{190,j}(250)$, averaged over 10000 \textcolor{black}{Monte Carlo} simulations starting from a fully connected network with $N=20$ nodes and for $\beta = 1, 0.7, 0.5$, as we move from left to right. \textcolor{black}{The Monte Carlo simulations were performed using an event-driven algorithm taking non-Markovianity into account.}
	The solid curves are the theoretical predictions as dictated by equation \eqref{eq:anal:istvan}.}
	\label{fig0}
	\end{figure}

	The probabilities 
	involving the counting process $N_{\beta}(t)$ have an explicit integral representation \cite{politi2011full} and for numerical purposes they 
	can be approximated well either with Monte Carlo simulations 
	or with a numerical integration scheme. 
	Once the transition probabilities of the embedded Markov 
	chain are known, every term is known in \eqref{eq:anal:istvan} and it can be used to exactly calculate the non-equilibrium probabilities $p_{i,j}(t)$ (\textcolor{black}{see Appendix B}).
	The excellent agreement between theory and simulation is shown in Fig. \ref{fig0}. 
	
	An immediate application is to use equation \eqref{eq:anal:istvan} and compute theoretically and numerically 
	the expected number of active links in the network at a given time. Starting from any initial number of active links 
	$i_0$, use \eqref{eq:anal:istvan} to compute 
	\begin{align}\label{eq:expQ}
	\E^{i_0}(X(t)) &= \bE^{i_0}\! \big(\textbf{e}_{i_0}^T Q^{N_{\beta}(t)}\textbf{v}_{0,M}\big)\notag \\
	&= \textbf{e}_{i_0}^T\bE^{i_0}\!\big(  Q^{N_{\beta}(t)}\big)\textbf{v}_{0,M}.
	\end{align}
	In the equation above 
	$\textbf{e}_{i_0}$ is the standard basis factor with $1$ at the $i_0$-th coordinate and  
	$\textbf{v}_{0,M} = (0,1,2, \ldots, M)^T$. Note that in the particular case where $Q$ diagonalizes, 
	the analytical expression for the expectation \eqref{eq:expQ} 
	is merely a linear combination of different values
	of the probability generating function of $N_{\beta}(t)$, $G^{(\beta)}(s;t)$, given by (see \cite{laskin2003fpp})
	$
		G^{(\beta)}(z;t) = \bE(z^{N_{\beta}(t)}) = E_{\beta}((z-1)t^{\beta}).
	$
	Note that when $Q$ diagonalises, there is no need for 
	simulating a large number of realisations to estimate the expectation; a fast numerical integration scheme 
	is sufficient.
	
%	This exactly solvable model of non-Markovian network dynamics is an important, tractable special case 
%	 of a more general theory for non-Markovian 
%	processes outlined in \cite{haenggi1981old} and 
%	it is related to recent outstanding developments in probability theory. 
%	In fact, up to this point, we provided a bottom-up derivation for the master equation of some 
%	fractional birth and death processes 
%	in a finite capacity system, 
%	introduced in \cite{orsingherpolito}.
	%\textbf{What did we learn from this part, how are we going to use it and why it is important??? We need to say something her to conclude!!!}

\subsection{Approximation of Pareto-distributed inter-event times}
	
 We now compare the behaviour between two RLAD networks; one with 
		Mittag-Leffler times and one where we 
	alter the waiting time distribution to a generalised Pareto($\delta$) with density
	\be
	f_T(t) = \frac{\delta - 1}{(1 + t)^{\delta}}, \quad t >0.
	\ee 
	Exponent $\delta = 1 + \beta$ in order for the tails of Mittag-Leffler and the 
	Pareto to have the same behaviour at infinity. 
	In fact we compare the two networks over three layers of increasing complexity.
	First, in Fig.~\ref{fig:1}(b,c) we plot the probability mass function for the number of singly-counted links 
	at a pre-specified time horizon $T=2000$, averaged over 
	5000 simulations. As a point of reference, output from the Markovian RLAD is shown in Fig.~\ref{fig:1}(a).
	The theoretical curve at equilibrium is the large time limit of equation \eqref{eq:anal:istvan} and it is the mass function of a binomial distibution
	with $M$ trials and success probability $1/2$. Second, 
	in Fig.~\ref{fig:1}(d) we plot $\E(X(t))$ as a function of time up to time 
	$T=2000$. The two curves could also be computed based on equation \eqref{eq:expQ}. 
	
	The excellent agreement is
	achieved by finding a suitable scaling $\gamma$ so that the c.c.d.f.~(survival functions) 
	of the two distributions match well, at least up to the pre-specified time horizon that we want to study. 
	Further details can be found in \textcolor{black}{Appendix C}. The matching is good for 
	$\beta < 0.9$ while for larger $\beta$, this idea can be used to 
	study the stochastic dominance between the two coupled networks and offer rigorous bounds. 

\subsection{Markovian SIS on non-Markovian RLAD}
	
	Finally, we compare the two network \textcolor{black}{dynamics} indirectly, when we allow a Markovian epidemic to run 
	while the networks evolve.  As discussed above, human activity tends to be bursty and non-Markovian \cite{barabasi2010bursts}. 
	During an epidemic, individuals become wary of the risk posed by it and
	one way to avoid infection is by limiting or reducing their number of contacts.
	This justifies the deletion of links as time evolves. On the other hand, close contacts 
	cannot realistically be removed and some level of communication and social cohesion must be maintained.
	Such behaviour in activation-deletion is not necessarily Markovian in nature, thus alternative 
	non-Markovian dynamic network models are necessary. 

%	Perhaps a more revealing motivating example would be 
%	the propagation of a `viral clip'  on online networks, (for example Facebook, since it seems 
%	the average graph distance between two users is less than four). 
%	links  correspond to friends and a deletion of a link 
%	can be considered extremely slow and rare during the time frame that the clip survives.
%	
	
%	While the present network model is simple, it does 
%	capture and accounts for some of the main features described above.
	
%	Hence, we combine a Markovian SIS model with the 
%	non-Markovian RLAD network model. 	
%	and this will have important consequences for how the re-shaping of contacts in the case of an epidemic
%	will impact on the severity of the outbreak. The epidemic starts spreading but individuals become wary of the infection so 
%	often they try to avoid it if they are healthy (or protect others if they are not). 
%	This justifies the deletion of links as time evolves. 
%	Often however one is obliged to communicate 
%	with others so links might  be activated again. Such behaviour in activation-deletion 
%	is not necessarily Markovian in nature, 
%	thus other models with analytical tractability are necessary. 

	Nodes in the network represent individuals from 
	a population and links describe the contact patterns amongst these. 
	Each individual can be either infected ($I$) os susceptible ($S$). 
%	An infected individual can be healed at rate $\tau_H$ while a susceptible individual 
%	in contact with an infected one gets infected at rate $\tau_I$.
\textcolor{black}{An infected individual remains infected for exponentially distributed periods of time $T_H$ i.e. $T_H \sim \mathrm{Exp}(1/\tau_H)$, where
$\tau_H$ is the average time in which infectious individuals are healed. Similarly, infection occurs at the points of a Poisson process
with time to infection $T_I$ exponentially distributed, i.e. $T_I \sim \mathrm{Exp}(1/\tau_I)$, where $\tau_I$ is the average time in
which an infection spreads across a link connecting a susceptible and an infected node. In this framework, both network and epidemic dynamics can be considered in the context of event-driven simulations, where the timing of the next state change is always determined by the smallest waiting time and the precise event corresponding to it.}
The epidemic does not interfere with the network dynamics, however its propagation 
	is intertwined with the background dynamic network topology.  
	Initially, before the infection starts spreading, we assume that all links are present, 
	in order to avoid early stochastic extinction.
%	This is not problematic, 
%	indeed any initial condition is theoretically valid; 
%	however starting from a sparse network and adding links slowly will 
%	kill the epidemic fast. Naturally a rigorous study of all these effects are of interest, 
%	but escape the scope of the current  letter.  
	
\begin{figure*}
%\begin{center}

\hfill \includegraphics[height=4cm]{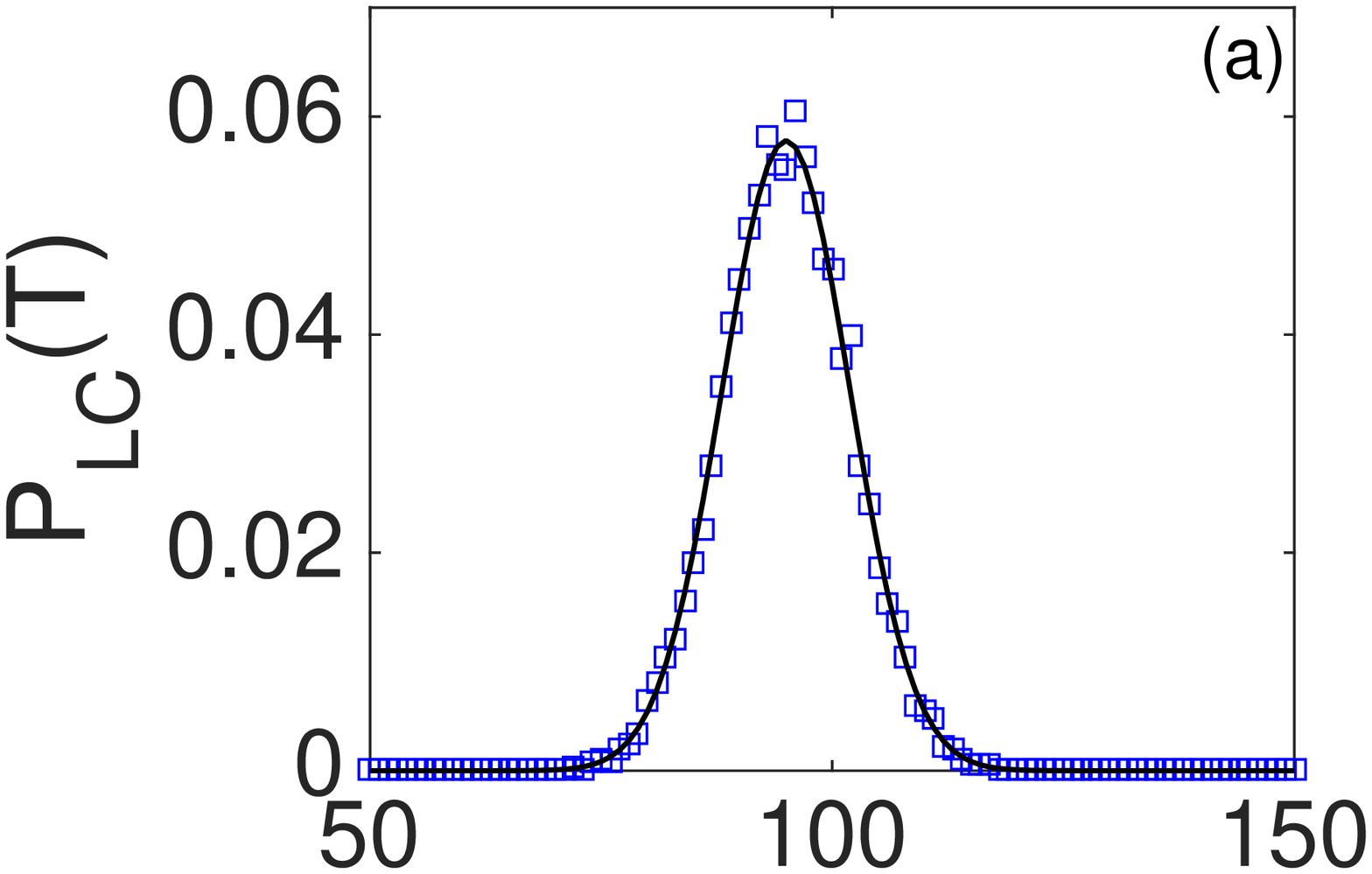} \hfill  \includegraphics[height=4cm]{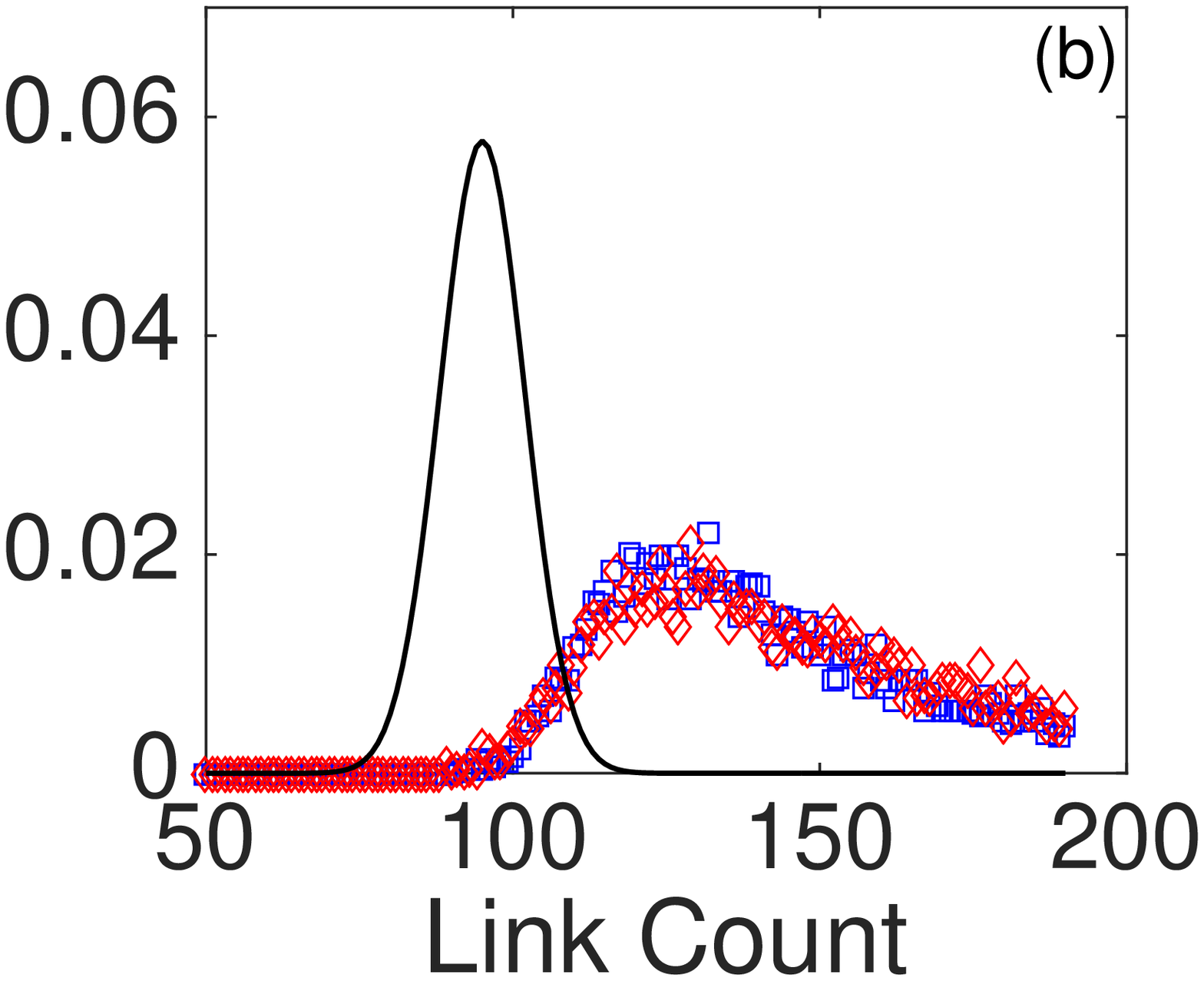} \hfill
  \includegraphics[height=4cm]{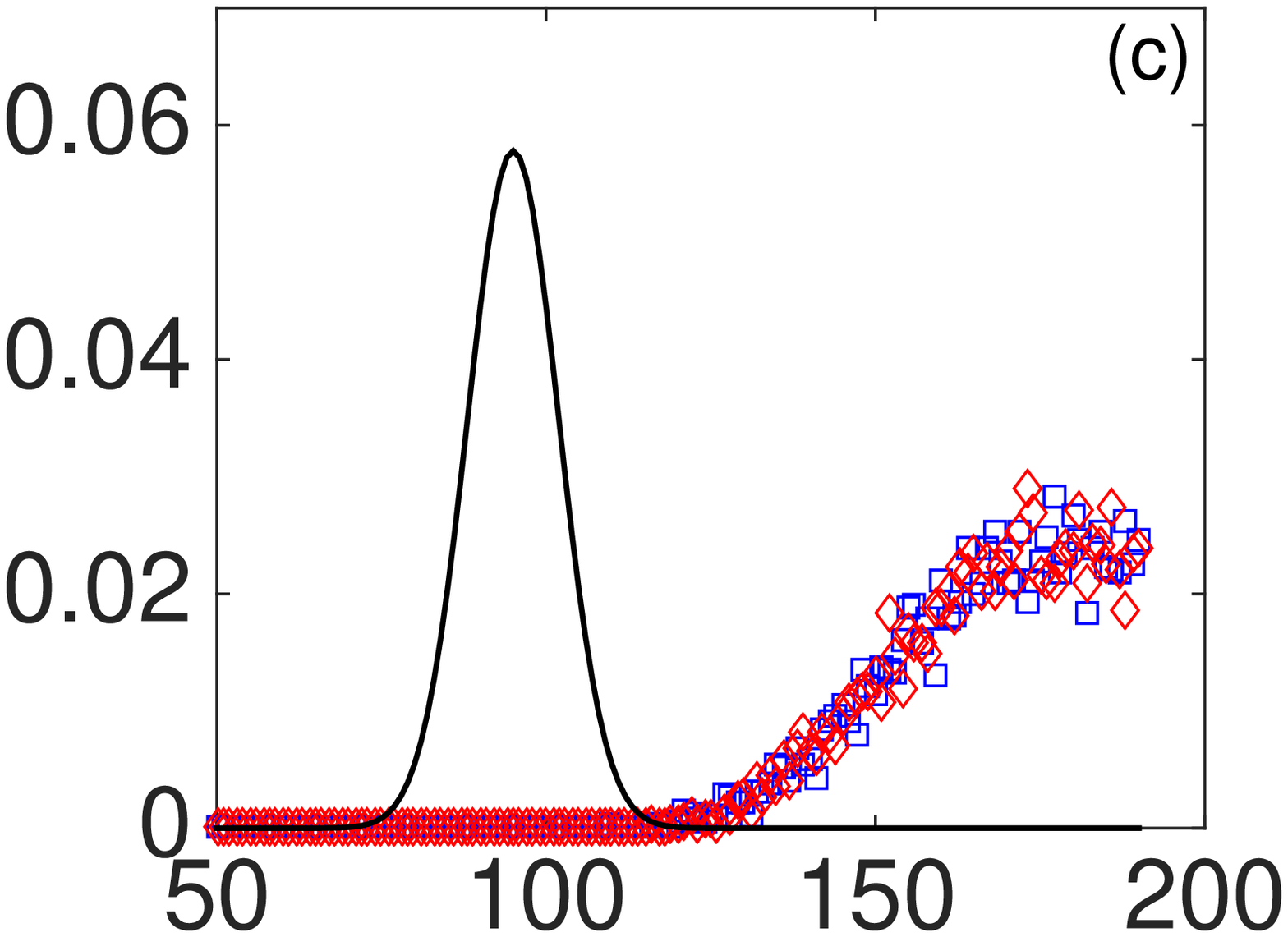}
\vfill
 \includegraphics[width=7cm, height=4.9cm]{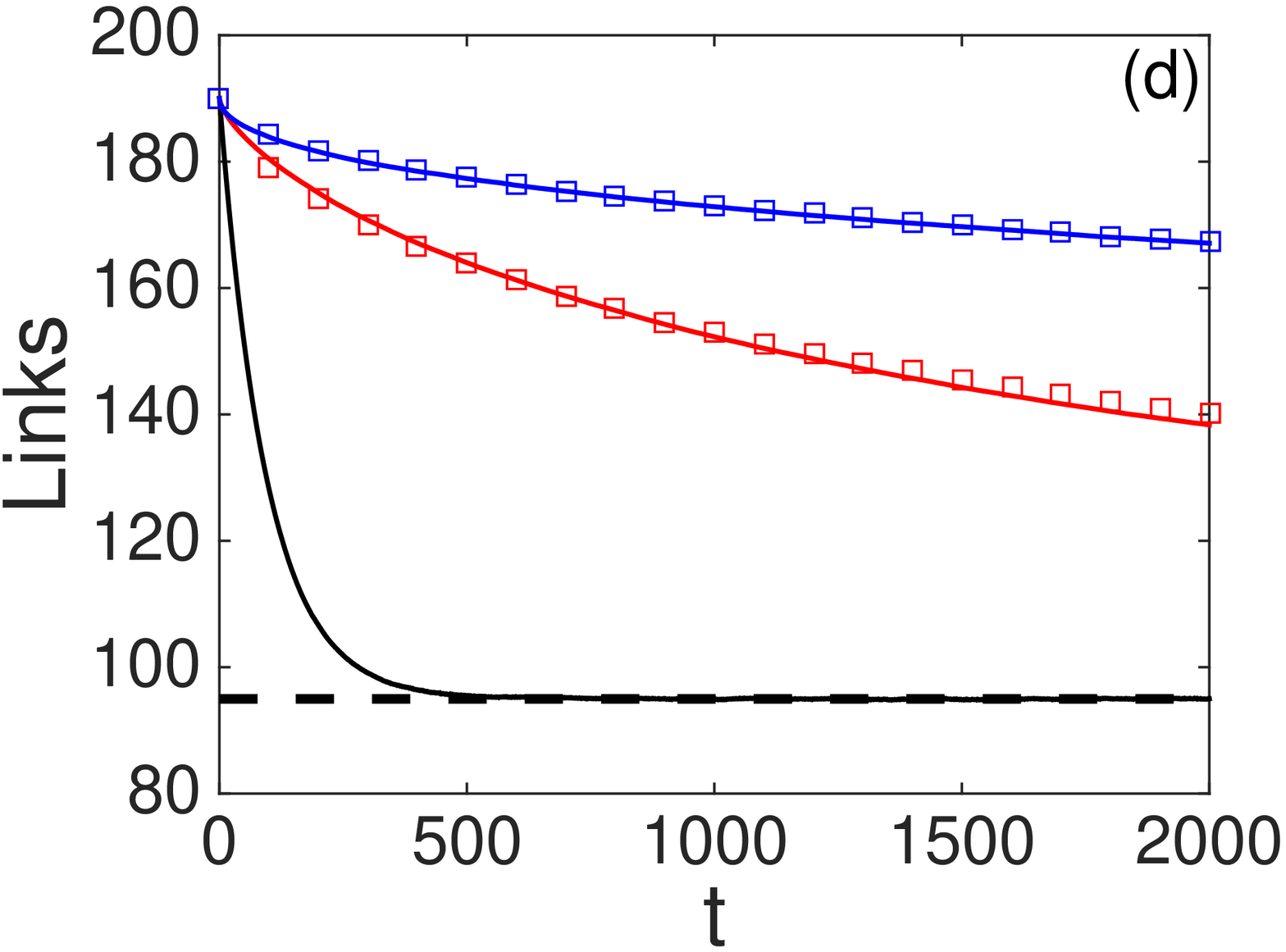} \quad\quad 
 \includegraphics[width=7cm, height=4.9cm]{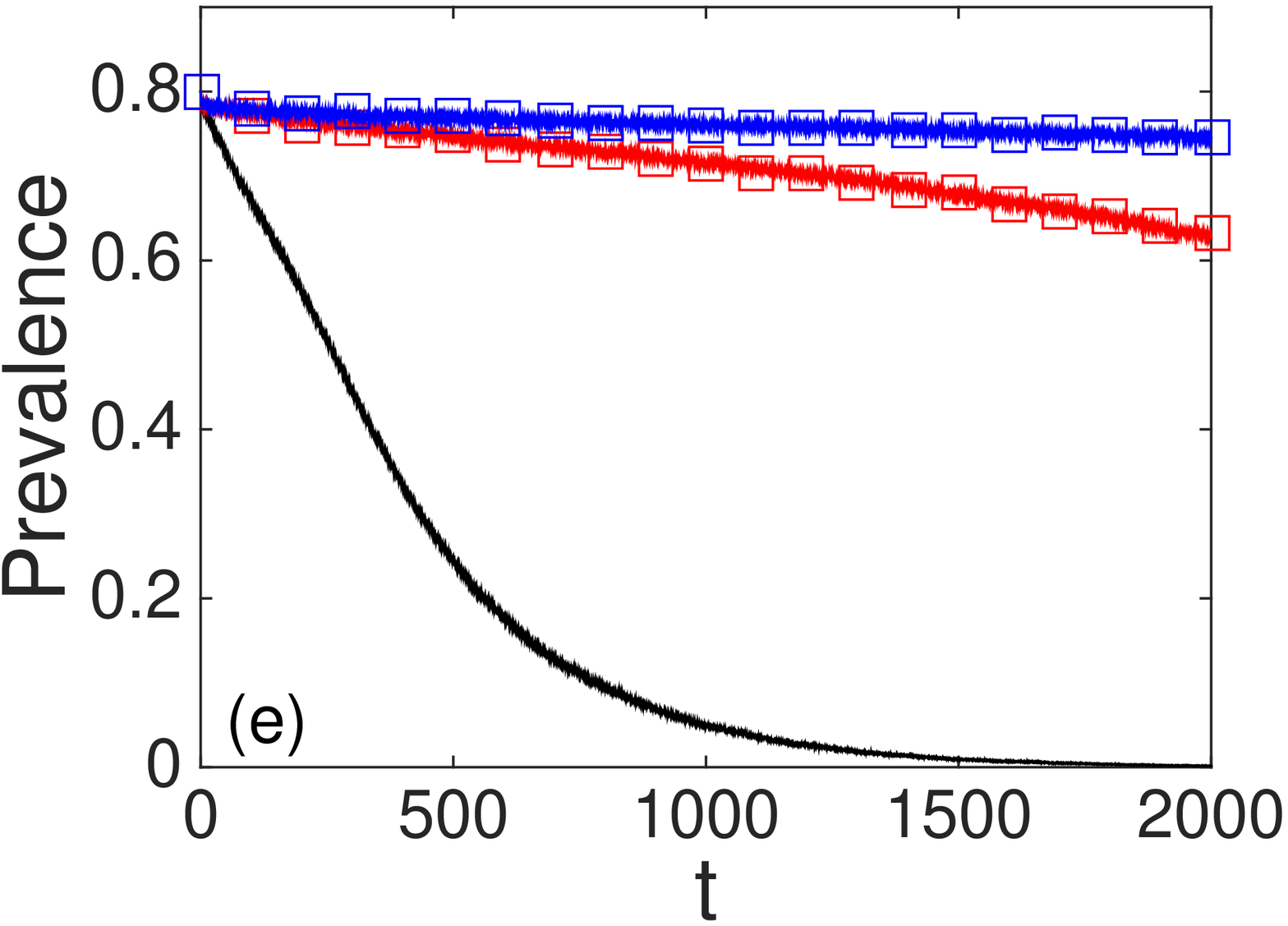} 

 \caption{(Color online) \textcolor{black}{Top: Distribution of links at $t=2000$ and comparison between Mittag-Leffler and Pareto.} Panels (a),(b),(c) show the distribution of singly counted link numbers at $t=2000$, 
 with theoretical equilibrium prediction (continuous line) and simulation results (discrete markers). \textcolor{black}{Theoretical equilibrium is the same
for all values of $\beta \in (0,1]$ and is given by the binomial distribution \eqref{eq:inv} here drawn as a continuous line for ease of representation.}
 $\square$ markers in the figures (a), (b), (c) correspond to $(\beta, \gamma)=(1,1),( 0.7,3.14), (0.5,4),$ respectively and $\lozenge$ markers for the corresponding Pareto($\delta$). 
 \textcolor{black} {Bottom:} On the second row, we show 
 (d) the expected number of singly counted links in the network (the dashed line is the theoretical prediction of $N(N-1)/4 = 95$ and (e) the expected prevalence.
 The continuous (noisy) line is for $(\beta, \gamma)=(1,1),( 0.7,3.14), (0.5,4),$ respectively from bottom to top and $\square$ markers are the corresponding values for the matched Pareto network. 
 The networks have $N=20$ nodes, and simulations start with a completely connected network. The $SIS$ epidemic is simulated as a Markovian process with 
 transmission and recovery rate $\tau_I=0.25$ and $\tau_H=1$, respectively. The spreading process initialises with $5$ infectious nodes. \textcolor{black}{Since we start with a fully connected network and a slow network dynamics, the prevalence rapidly increases from 25$\%$ to almost 80$\%$. This also reflects the relation between the time scales of the network and epidemic dynamics, with the epidemics being much faster in this example.}
The simulation is event-driven and is implemented by keeping track of all events and their waiting times. These are averaged over 5000 simulations and use 
$\alpha = 0$ (periodic case),  so all events create a change in the network \textcolor{black}{and equilibrium is reached earlier}.}
\label{fig:1}
%\end{center}
 \end{figure*}

	The simulations in Fig.~\ref{fig:1} show the prevalence (proportion of infected individuals (Fig.~\ref{fig:1}(e)) 
	on a Mittag-Leffler$_\gamma(\beta)$ RLAD network (solid lines) and a direct comparison (square markers) 
	with the Pareto($\delta$). Again, we use the same sets of $\beta, \gamma, \delta$ as before and 
	we emphasise the excellent agreement between the two. 
%	\be
%	f_T(t) = \frac{\alpha - 1}{(1 + t)^{\alpha}}, \quad t >0.
%	\ee 
%	Exponent $\alpha = 1 + \beta$ in order for the tails of Mittag-Leffler and the 
%	Pareto to have the same behaviour at infinity. 
%	The excellent agreement of the prevalence for the two epidemics is 
%	achieved by finding a suitable scaling $\gamma$ so that the c.c.d.f  (survival functions) 
%	of the two distributions match well, at least up to the pre-specified time horizon that we want to study. 
%	Further details can be found in the Supplemental material. This matching works well for 
%	$\beta < 0.9$ while for larger $\beta$ this idea can be used to 
%	study the stochastic dominance between the two coupled networks and to offer rigorous bounds. 
	
%	For the same sets of $\alpha, \beta, \gamma$ used in the epidemics, 
%	we further compare and note the agreement between 
%	the Mittag-Leffler and the Pareto RLAD networks 
%	by averaging the number of links for both networks up to time $t = 2000$ 
%	over 5000 simulations (Fig. 2b). The two curves are theoretically 
%	given by equation \eqref{eq:expQ}. 
%	
%	The expected value of the process 
%	is not the only statistic of the two networks that has agreement, as it can be seen from 
%	Figs.~(2c, 2d) where we plot the simulated probability densities  
%	for both networks when $\beta = 0.7$ and $0.5$ respectively.           
	
	Incidentally, as expected, the non-Markovian network dynamics create a striking effect by 
	slowing down the network dynamics and thus effectively 
	blocking the attainment of statistical equilibrium in a realistic time horizon (Fig.~\ref{fig0}, \ref{fig:1}(a,b,c)). 
	This leads to a heightened level of infectiousness in the population (Fig.~\ref{fig:1}(e))  and 
	highlights the importance of quick reactions. Naturally the statistical equilibrium will 
	be reached after a much longer time period, but the delayed curves 
	can now  be theoretically computed or approximated.
%	Notice the excellent agreement in the distribution of
%	the number of links in the network for the Markovian case (Fig. 1c) after a relatively short simulation 
%	whereas this does not happen in the non-Markovian case (Fig.1d,1e). 
%	Results on the convergence to equilibrium of these stochastic processes follow.
%\textit{Can we refer back to any formulas to claim agreement???}
	%The state of the epidemic process is binary 
	%vector where a $1$ at position $i$ at time $t$ indicates infection of the $i$-th individual. 
	%One is interested in the number of infected individual at a given time $t$. 
	%When the population graph is fixed, these dynamics are Markovian and time homogeneous.
	%The process has an absorbing state of all $0$'s, i.e. the disease died out 
	%and the entire population is just susceptible.
%	In our current setting we introduce a non-Markovian 
%	time inhomogeneous dynamic network. On a heuristic level, the choice of the joint 
%	dynamics can be justified is as follows.
%	The epidemic starts spreading but individuals become wary of the infection so 
%	often they try to avoid it if they are healthy (or protect others if they are not). 
%	This justifies the deletion of links as time evolves. 
%	Often however one is obliged to communicate 
%	with others so links might  be activated again. Such behaviour in activation-deletion 
%	is not necessarily Markovian in nature, 
%	thus other models with analytical tractability are necessary. 
%\section{Comparison to the Markovian RLAD model.} 
One way to explain this delayed convergence to equilibrium is to look at the mixing time of the embedded chain in the total variation distance. 
%This idea was considered  for some non-Markovian dynamic networks, for example in \cite{S-W-P-2015}.
%As evidenced by the simulations, in the non-Markovian case it is incredibly difficult for a large system 
%to reach the steady state in a sensible time frame. 
%This is because the Mittag-Leffler distribution is heavy tailed and 
%has infinite expectation for all parameter values $\beta \in (0,1)$. 
%Here, the expected number of events up to time $t$ with Mittag-Leffler($\beta$) is \cite{laskin2003fpp} $\E(N_{\beta}(t)) = O(t^{\beta})$. 
%On the other hand, using the law of large numbers, Theorem 1.1 and Example 4.3 
%in \cite{Che-Sal-13}, the embedded chain is mixed and close to equilibrium after $O(M^2\log M)$ time events. For the Mittag-Leffler RLAD, this number of events only happens on average after $O(M^{2/\beta}(\log M)^{1/\beta})$ time elapses.
%On the other hand, consider the Markovian RLAD with mean $1$, 
%exponentially distributed waiting times.  The equations are the same as \eqref{eq:final0}, \eqref{eq:final1}, \eqref{eq:final2} with the usual, rather than the Caputo derivative. 
To be more specific, the number of active links is a continuous time, irreducible birth-death chain with a unique binomial invariant 
distribution $\pi$, independent of the delay parameter $\alpha$, given by 
	\be\label{eq:inv}
		\pi_k= \lim_{t \to \infty}\P\{ X(t) = k|X(0)=i\} = { M \choose k}2^{-M}, \quad k \in \cS. 
	\ee
\textcolor{black}{This can also be deduced from the fact that in the aperiodic $\alpha$-delayed case, at equilibrium, individual graphs are uniformly distributed. That is because the chain on the set of distinct graphs has a doubly stochastic transition matrix. With this in mind, the degree distribution of a single node in network chosen uniformly at random can be immediately computed as follows. Let $v_1$ be a selected 
node in $G$ and define $G_{v_1}$ the subgraph of $G$ where $v_1$ and all its incident links are deleted. $G_{v_1}$ is now a graph on the set $\{v_2, \ldots, v_N\}$ that has at most $K = {N-1 \choose 2} = M - N+1$ links. Let $ \text{deg}(v_1)$ denote the degree of $v_1$. Then,
	\begin{align}
		\P\{ \text{deg}(v_1) = \ell \} 
			&= \!\!\!\!\sum_{\text{graphs } G: \text{ deg }(v_1) = \ell  } \!\!\!\!\!\!\!2^{-M}\notag\\
			&= 2^{-M} \text{card}\{ \text{graphs } G: \text{ deg }(v_1) = \ell \}\notag\\
			&=2^{-M}{ N-1 \choose \ell} \text{card}\{ \text{graphs }G_{v_1} \}\notag \\
			&=2^{-K-N+1}{ N-1 \choose \ell}\sum_{i=0}^K{K \choose i}\notag\\ 
			&= { N-1 \choose \ell}2^{-N+1}. \label{eq:deg:inv}
	\end{align}
	The third equality is a counting argument. The number of graphs such that $v_1$ has exactly $\ell$ incident links, is constructed  by first selecting where those links go, and then constructing the subgraph $G_{v_1}$.
This can be understood heuristically as follows. Select a link and focus on one of the two nodes. If this node has $h$ active links, this number will either go to $h+1$ with probability $(N-1-h)/(N-1)$ or to $h-1$ with probability $h/(N-1)$. This leads to an invariant binomial degree distribution \eqref{eq:deg:inv},
%\begin{equation}
%\mathbb{P} (\mathrm{degree}=l) = {N-1 \choose l} \frac{1}{2^{N-1}},
%\end{equation}
with an average degree of $(N-1)/2$, which amounts to $N(N-1)/4$ active edges in the network in line with the average of the link distribution from \eqref{eq:inv}.}
The chain mixing time $t_{\text{mix}}(\e)$ is the minimal time so that the total variation distance 
between the measures $\pi$ and $p(t)$ is smaller than some tolerance $\e$, i.e. 
	\be
		\| p(\tmix(\e))-\pi \|_{\text{TV}} = \sup_{k \in \cS}|p_{k}(\tmix(\e)) - \pi_k| < \e.
	\ee
For the Markovian RLAD with $\alpha > 0$, use Theorem 1.1 and Example 4.3 in \cite{Che-Sal-13}
%\textcolor{black}{(Look at diaconis 1988 monograph, exact computations)} 
to see 
	\be
		 \tmix(\e) \le C\e^{-2} M^2\log M.
	\ee
Thus, the Markov chain approximates relatively well its equilibrium by time 
$C\e^{-2} M^2\log M$. In particular this implies that, on average, the Markovian RLAD continuous chain 
needs $O(M^2\log M)$ time, and therefore by the law of large numbers, it needs this order many events until it is well mixed.
In fact this bound is also true for the embedded discrete chain. 
This $n = M^2\log M$ should be considered as a \emph{necessary} lower bound 
of steps so that the sample average of the probabilities of the embedded chain approximates 
$\pi$. Therefore in order to have an acceptable level of accuracy for the embedded chain 
when the RLAD has Mittag-Leffler waiting times, using \eqref{eq:m-l:number}, 
we need  a higher polynomial order 
$O(M^{2/\beta}(\log M)^{1/\beta})$ time in order to 
guarantee on average the same number of events, and thus to guarantee 
a near equilibrium behaviour for the embedded chain.  
Study of the slow-down phenomenon for non-Markovian dynamic networks, 
using the total variation distance, can be 
found for example in \cite{S-W-P-2015}.

\section{Discussion and conclusions}

\subsection{Generalization}

\textcolor{black}{Equation \eqref{eq:anal:istvan} can be generalized to any counting process and any discrete Markov chain and, as a consequence, to any embedded Markovian graph dynamics \cite{raberto2011graphs}. To be more specific, let
$q_{i,j}$ denote the one-step transition probability from state $i$ to state $j$ for a discrete Markov chain $X_n$ and let $N(t)$ be a generic counting renewal process. Then, for the process
\begin{equation}
\label{generalisation}
X(t) = X_{N (t)} = X_n 1\!\!1\{ S_n \le t < S_{n+1} \},
\end{equation}
the probabilities $p_{i,j} (t) = \mathbb{P} \{ X(t) =j | X(0) = i \}$ are given by
\begin{equation}
\label{generalisation1}
p_{i,j}(t) = \overline F_T(t)\delta_{ij} + \sum_{n=1}^\infty q^{(n)}_{i,j} \P\{N(t) = n\},
\end{equation}
where the symbols have the same meaning as in \eqref{eq:anal:istvan} and $\{T_i \}_{i=1}^\infty$ is a sequence of i.i.d. positive random variables with the usual meaning of inter-event times with arbitrary distribution, not necessarily with fat tails and infinite mean. The reader is invited to follow the first proof of Appendix B by replacing $N_\beta (t)$ with a generic counting renewal process. This will convince the reader of the wide generality of this result. A heuristic argument to justify \eqref{eq:anal:istvan} and \eqref{generalisation1} runs as follows. In the time interval $(0,t)$, $n \geq 0$ events may have occurred. In the case $n=0$, at time $t$ the process is still in state $i$ and $\mathbb{P}\{N(t) =0\} = \mathbb{P}\{T>t\} = \overline F_T(t)$. If $n \geq 1$, the probability of being in state $j$ after $n$ events is given by $q^{(n)}_{i,j}$. Given the independence between the renewal process and the Markov chain, the probability of being in state $j$ at time $t$ and $n$ transitions occurring in the time interval $(0,t)$ is $q^{(n)}_{i,j} \P\{N(t) = n\}$. Now, all these events are exhaustive and mutually exclusive. Then, total probability and infinite additivity imply that $p_{i,j} (t)$ is given by \eqref{generalisation1}. These considerations suggest further generalizations taking into account possible dependence within the couple $\{(X_n,T_n)\}_{n \geq 1}$ as well as serial dependence or state dependence of inter-event times.}

\subsection{Example: A simple probabilistic model for relaxation in dielectrics}

\textcolor{black}{In order to illustrate the generalization discussed above with an example, we consider relaxation phenomena \cite{weron}. Probabilistic modelling of relaxation assumes that a physical system (e.g. a molecule) can exist in two states $A$ and $B$. We further assume that state $A$ is transient and state $B$ absorbing, so that the deterministic embedded chain has the following transition probabilities $q_{A,A} = 0$, $q_{A,B} = 1$, $q_{B,A}=0$, and $q_{B,B} = 1$. This means that if the system is prepared in state $A$, it will jump to state $B$ at the first step and it will stay there forever. Now suppose that the inter-event time $T$ is random and follows an exponential distribution with rate $\lambda =1$ for the sake of simplicity. Based on equation \eqref{generalisation1}, we immediately have $p_{A,A} (t) = \overline F_T(t)=\exp(-t)$. Therefore, the probability of finding the system in the initial state decays exponentially towards zero. This relaxation function is the solution of
\begin{equation}
\frac{d}{dt} p_{A,A} (t) = - p_{A,A} (t), \; \; p_{A,A} (0) =1. 
\end{equation}
The response function is defined as $\xi_D (t) = -d p_{A,A} (t)/dt$ and its Laplace transform is $1/(1+s)$. For $s=-i \omega$ this is the
Debye model \cite{weron}. If inter-event times follow the Mittag-Leffler distribution, we get $p_{A,A} (t) = \overline F_T(t)=E_\beta (-t^\beta)$. This is the solution of \cite{scalas2004}
\begin{equation}
\frac{d^\beta}{dt^\beta} p_{A,A} (t) = - p_{A,A} (t), \; \; p_{A,A} (0) =1. 
\end{equation}
In this case, the Laplace transform of the response function $\xi_{CC} (t) = -d p_{A,A} (t)/dt$ is $1/(1+s^\beta)$ and for $s=-i \omega$, we get the Cole-Cole model \cite{colecole1,colecole2,weron}.}

\subsection{Final considerations}

In conclusion, we provide an exactly solvable non-Markovian dynamic network model.
The RLAD is particularly
attractive as it has analytical and numerical tractability coming from fractional calculus. 
We are able to explicitly use the master equation formalism and 
%to show that
%under appropriate transformations, the master equation can be made equivalent to 
%a problem posed in fractional calculus. We are able to 
analytically derive the distribution of the 
number of links in the network for arbitrary times $X(t)$, consequently computing $\E(X(t))$.
We highlight an important connection and possible avenue to approximate non-Markovian problems
using fractional calculus, by coupling a Pareto network and show the agreement with the tractable model.
\textcolor{black}{Moreover, we discuss how our result can be extended to a generic counting renewal process.}

\begin{acknowledgments}
\textcolor{black}{This paper was partially supported by an SDF start-up research grant provided by the University of Sussex.}
\end{acknowledgments}

\section*{Appendix}

In this Appendix, we cover the rigorous proofs of the equations shown in the main text and further 
clarify some notions. Some details about the procedure used to couple the Pareto distribution with the Mittag-Leffler are also highlighted.

\subsection*{A. Derivation of fractional equations.}

\noindent 
We want to show that equations \eqref{eq:final0}, \eqref{eq:final1} and \eqref{eq:final2} in the main text are obtained from \eqref{eq:master}.
%\begin{align} 
%		p_{i,j}(t) \!= \overline F^{(\beta)}_T(t)\delta_{ij} 
%			+\! \sum_{\ell \in \cS} q_{\ell, j} \!\! \int_0^t \!\! p_{i,\ell}(u)f^{(\beta)}_{T}(t-u)\,du,
%			\label{eq:master0}
%	\end{align}
%where $p_{i,j}(t) = \P\{X(t) = j | X(0) = i\}$, $\overline F^{(\beta)}_T(t) = 1-F^{(\beta)}_T(t)$ is the tail (complementary cumulative distribution function) and 
%$f^{(\beta)}_{T} (t)$ the Mittag-Leffler density or order $\beta$.
The analysis proceeds by way of Laplace transforms. They are defined as
\be
\cL (g(t);s) = \int_0^\infty dt \, g(t) \, \mathrm{e}^{-st}
\ee
for a suitable function $g(t)$.
In the case of the Mittag-Leffler distribution defined in the main text,
we have 
	\be\label{eq:ml-lap}
		\cL\left(  \overline F^{(\beta)}_T(t) ; s \right) = \frac{s^{\beta-1}}{1 + s^{\beta}}  \text{ and } 
		\cL\left( f^{(\beta)}_T(t) ; s \right) = \frac{1}{1 + s^{\beta}}.
	\ee
	
For the computation that follows we use the symbol $\tilde g(s)$ 
to denote the Laplace transform $\cL(g;s)$ of any function $g$. Taking the Laplace 
transform of \eqref{eq:master} and using equations (4), (5), (6) in the main text for our particular example, we have 
for $1 \le j \le M-1$
	\begin{align} \label{eq:Lmaster:bulk}
		\tilde p_{i,j}(s) &= \widetilde{ \overline F}^{(\beta)}_T(s)\delta_{ij} + \tilde f^{(\beta)}_{T}(s)\alpha \tilde p_{i,j}(s)  \\
			&+ \tilde f^{(\beta)}_{T}(s)(1-\alpha)\notag \\
			&\times \Big( \frac{M -j +1}{M} \tilde p_{i, j-1}(s) +  \frac{j+1}{M}\tilde p_{i,j+1}(s) \Big).\notag
	\end{align}
The boundary cases $j =0, j=M$ have Laplace transforms 
	\begin{equation}
		\tilde p_{i,0}(s) = \widetilde{ \overline F}^{(\beta)}_T(s)\delta_{i0} 
				+  \tilde f^{(\beta)}_{T}(s)\left(\frac{1-\alpha}{M}\tilde p_{i,1}(s) + \alpha \tilde p_{i,0}(s)\right),
			\label{eq:Lmaster:j0}
	\end{equation}
	\begin{align}
		\tilde p_{i,M}(s) &= \widetilde{ \overline F}^{(\beta)}_T(s)\delta_{iM} \notag \\
				&+  \tilde f^{(\beta)}_{T}(s)\left(\frac{1-\alpha}{M}\tilde p_{i,M-1}(s) + \alpha \tilde p_{i,M}(s)\right)
			\label{eq:Lmaster:jM}
	\end{align}
	respectively. We finish the computation starting from \eqref{eq:Lmaster:bulk}, in the case where $ 1 \le j \le n-1$. 
	The remaining  cases 
	follow similarly. Multiply both sides of \eqref{eq:Lmaster:bulk} by $s$ and then 
	subtract $p_{i,j}(0)=\delta_{ij}$ from both sides. Then, using \eqref{eq:ml-lap}, equation \eqref{eq:Lmaster:bulk} becomes 
	\begin{align*}
		\cL\Big(& \frac{d p_{i,j}(t)}{dt}; s \Big) = s\widetilde{ \overline F}^{(\beta)}_T(s)\delta_{ij} -p_{i,j}(0)
		+ \frac{s}{1+ s^{\beta}}\alpha \tilde p_{i,j}(s)\\
			&\phantom{x}
				+ \frac{s(1-\alpha)}{1+ s^{\beta}}\Big( \frac{M -j +1}{M} \tilde p_{i, j-1}(s) + \frac{j+1}{M}\tilde p_{i,j+1}(s) \Big)\\
			&= \frac{s^{\beta}}{1+s^{\beta}}\delta_{ij} -\delta_{ij} + \frac{s}{1+ s^{\beta}}\alpha \tilde p_{i,j}(s)\\
			&\phantom{x}
				+ \frac{s(1-\alpha)}{1+ s^{\beta}}\Big( \frac{M -j +1}{M} \tilde p_{i, j-1}(s) + \frac{j+1}{M}\tilde p_{i,j+1}(s) \Big),
	\end{align*}
	thus, after some algebraic manipulations we have %the equivalent equation
	\begin{align} \label{eq:almost} 
		&\frac{1+s^{\beta}}{s}\cL\Big( \frac{d p_{i,j}(t)}{dt}; s \Big) = \frac{-\delta_{ij}}{s} + \alpha \tilde p_{i,j}(s)\\
		&\phantom{x}+ (1-\alpha)\!\left(\frac{M -j +1}{M} \tilde p_{i, j-1}(s) + \frac{j+1}{M}\tilde p_{i,j+1}(s)\right).\notag
	\end{align}
	Focus for the moment on the factor $s^{-1}(1 +s^{\beta})$. Its inverse Laplace transform is 
	\be\label{eq:caputo}
		\cL^{-1}\big( s^{-1}(1 +s^{\beta}); t \big) = \frac{t^{-\beta}}{\Gamma(1-\beta)} + 1 = \Phi_{\beta}(t) + 1.
	\ee
	Kernel $\Phi_{\beta} (t)$ is what is used in fractional calculus to define the \emph{ $\beta$ fractional Caputo derivative} (see reference [19] in the main text)	
	of a function $f(t)$, given by 
	\[
		\frac{d^{\beta} f(t)}{ d\, t ^{\beta}} = \int_{0}^t \Phi_{\beta}(t-t')\frac{d\,f(t')}{dt'}\,dt'. 
	\]
	Thus, use \eqref{eq:caputo} to write 
	the left hand side of \eqref{eq:almost} as a product of two Laplace transforms.
Then take the Laplace inverse of \eqref{eq:almost} to conclude 
	\begin{align}\label{al:final0}
		&\frac{d^{\beta} p_{i,j}(t)}{ d\, t ^{\beta}}  = -(1-\alpha)\ p_{i,j}(t) \\
			&\quad+ (1-\alpha)\left(\frac{M -j +1}{M} p_{i, j-1}(t) + \frac{j+1}{M}p_{i,j+1}(t)\right).\notag
	\end{align}
	Similarly, the equations of the boundary terms are derived \eqref{eq:final1}, \eqref{eq:final2}.
	\bigskip

\subsection*{B. Solution to the fractional equations.} 

The solution to equations (12), (13), (14) can be seen to be equation (15) in two different ways. 
One is the standard law of total probability, where the space is partitioned according to the number of jumps of the counting process $N_{\beta}(t)$:
\begin{align*}
	p_{i,j}&(t) = \P\{X(t) = j | X(0) = i\} \\
	&= \sum_{k=0}^\infty \P\{X(t) = j, N_\beta(t) =k | X(0) = i\}\\
	&= \sum_{k=0}^\infty \P\{X(t) = j | N_\beta(t) =k , X(0) = i\} \P\{ N_{\beta}(t) = k\}\\
	&= \P\{X(t) = j | N_\beta(t) = 0 , X(0) = i\} \P\{ N_{\beta}(t) = 0\} \\
	&\phantom{x}+ \sum_{k=1}^\infty \P\{X(t) = j | N_\beta(t) =k , X(0) = i\} \P\{ N_{\beta}(t) = k\}\\
	&= \P\{X(t) = j | T_1\ge t , X(0) = i\} \P\{ T_1 \ge t \} \\
	&\phantom{x}+ \sum_{k=1}^\infty \P\{X(t) = j | N_\beta(t) =k , X(0) = i\} \P\{ N_{\beta}(t) = k\}\\
% &= \delta_{ij} \P\{ T_1 \ge t \} \\
%	&\phantom{x}+ \sum_{k=1}^\infty \P\{X_k1\!\!1\{ S_k \le t < S_{k+1}\} = j |  1\!\!1\{ S_k \le t < S_{k+1}\} , X_0 = i\} \P\{ N_{\beta}(t) = k\}\\
	&= \delta_{ij} \overline F^{(\beta)}_T(t)\\
	&+ \sum_{k=1}^\infty \P\{X_k = j |  1\!\!1\{ S_k \le t < S_{k+1}\} , X_0 = i\} \\
	&\phantom{xxxxx}\times \P\{ N_{\beta}(t) = k\}\\
	&= \delta_{ij} \overline F^{(\beta)}_T(t)+ \sum_{k=1}^\infty \P\{X_k = j |  X_0 = i\} \P\{ N_{\beta}(t) = k\},
\end{align*} 
where it finally leads to 
\be\label{eq:s1}
p_{i,j}(t) = \delta_{ij} \overline F^{(\beta)}_T(t)+ \sum_{k=1}^\infty q^{(k)}_{ij} \P\{ N_{\beta}(t) = k\}.
\ee

Equation \eqref{eq:s1} is equation \eqref{eq:anal:istvan} in the main text and as we say, 
gives the theoretical solution to the fractional equations, because  of an explicit integral representation of  
$\P\{ N_{\beta}(t) = n\}$. It is given by 
\begin{align*}
	\P\{ N_{\beta}(t) = n\} &= \frac{t^{\beta n}}{n!}E_{\beta}^{(n)}(-t^\beta) \\ 
	&= \int_{0}^\infty F_{S_{\beta}}(t;u)\left( 1 - \frac{u}{n}\right) \frac{u^{n-1}}{(n-1)!}e^{-u}\,du.
\end{align*}
Function $F_{S_{\beta}}(t;u)$ is the c.d.f.~of a stable random variable $S_{\beta}(\nu, \gamma, \delta)$ with index $\beta$, skewness parameter $\nu=1$, scale $\gamma = (u\cos(\pi \beta/2))^{1/\beta}$ and location $\delta=0$. This integral representation was used to numerically compute the solid curve in Figure \ref{fig0} \cite{politi2011full}.

We now verify via Laplace transforms that this solution \eqref{eq:s1} indeed verifies the fractional equations. For simplicity we set the delay parameter $\alpha=0$ and we only show it for equation \eqref{eq:final0}. We need 
\[ 
\cL\left( \frac{d^{\beta} g}{dt^\beta}; s\right) = s^{\beta}\tilde g(s) - s^{\beta-1} g(0+), \,\,\]
\[  \cL\left( \P\{N_{\beta}(t) = n\}; s\right) =\widetilde{\overline F}^{(\beta)}_T(s) \left(\tilde f^{(\beta)}_{T}(s)\right)^n\!\!= \frac{\widetilde{\overline F}^{(\beta)}_T(s)}{(1+s^{\beta})^n}. 
\] 
The Laplace transform of \eqref{eq:final0} 
\begin{align*}
s^{\beta}&\tilde p_{i,j}(s) - s^{\beta-1} \delta_{ij}  \\
&= -\tilde p_{i,j}(s) + \frac{M -j +1}{M} \tilde p_{i, j-1}(s) + \frac{j+1}{M}\tilde p_{i,j+1}(s),
\end{align*}
or after an algebraic manipulation 
\be\label{eq:magic}
(1+s^{\beta}) \tilde p_{i,j}(s) =  s^{\beta-1} \delta_{ij} + q_{j-1,j} \tilde p_{i, j-1}(s) + q_{j+1,j}\tilde p_{i,j+1}(s).
\ee
To verify \eqref{eq:magic}, directly take the Laplace transform in \eqref{eq:s1} to write 
\be
\tilde p_{i,j}(s) = \delta_{ij} \widetilde{\overline F}^{(\beta)}_T(s)+ \widetilde{\overline F}^{(\beta)}_T(s)\sum_{k=1}^\infty q^{(k)}_{ij} \left(\tilde f^{(\beta)}_{T}(s)\right)^k
\ee
and substitute to the right hand side in \eqref{eq:magic} that now reads
\begin{align*}
 s^{\beta-1}&\delta_{ij} + q_{j-1,j} \tilde p_{i, j-1}(s) + q_{j+1,j}\tilde p_{i,j+1}(s) \\
 &= s^{\beta-1} \delta_{ij} \\
 &\,\,+ q_{j-1,j}\widetilde{\overline F}^{(\beta)}_T(s) \left(  \delta_{i,j-1} + \sum_{k=1}^\infty q^{(k)}_{i,j-1} \left(\tilde f^{(\beta)}_{T}(s)\right)^k \right) \\
 &\phantom{x}+ q_{j+1,j}\widetilde{\overline F}^{(\beta)}_T(s)\left(\delta_{i,j+1} + \sum_{k=1}^\infty q^{(k)}_{i,j+1} \left(\tilde f^{(\beta)}_{T}(s)\right)^k\right)\\
 &= s ^{\beta-1}\delta_{ij} + (q_{j-1,j}  \delta_{i,j-1} + q_{j+1,j}\delta_{i,j+1} )\widetilde{\overline F}^{(\beta)}_T(s)\\
 &+ \widetilde{\overline F}^{(\beta)}_T\!\!(s)\sum_{k=1}^\infty (q_{j-1,j}q^{(k)}_{i,j-1} +  q_{j+1,j}q^{(k)}_{i,j+1} )\left(\tilde f^{(\beta)}_{T}(s)\right)^k \\
  &= s^{\beta-1} \delta_{ij} + q_{i,j}\widetilde{\overline F}^{(\beta)}_T(s) \\
  &\quad+(1+s^{\beta}) \widetilde{\overline F}^{(\beta)}_T(s)\sum_{k=1}^\infty q^{(k+1)}_{i,j}\left(\tilde f^{(\beta)}_{T}(s)\right)^{k+1} \\
   &= s^{\beta-1} \delta_{ij} + q_{i,j}s^{\beta-1} \tilde f^{(\beta)}_{T}(s)\\
   &\quad+(1+s^{\beta}) \widetilde{\overline F}^{(\beta)}_T(s)\sum_{k=1}^\infty q^{(k+1)}_{i,j}\left(\tilde f^{(\beta)}_{T}(s)\right)^{k+1} \\
   &= (1+s^{\beta}) \widetilde{\overline F}^{(\beta)}_T(s) \left[\delta_{ij} +
   +\sum_{k=1}^\infty q^{(k)}_{i,j}\left(\tilde f^{(\beta)}_{T}(s)\right)^{k}\right] \\
   &=(1+s^{\beta}) \tilde p_{i,j}(s), 
\end{align*}
which is the left hand side of \eqref{eq:magic}.

\subsection*{C. Stochastic coupling with the Pareto distribution.}

Now we show how the complementary cumulative distribution functions of the Pareto($\delta$) distribution and the Mittag-Leffler$_\gamma(\beta)$ of the same exponent can be shown to match just by manipulating the scaling factor $\gamma$.  

\begin{figure}[h]\begin{center}
	\includegraphics[width=9cm,height=5.75cm]{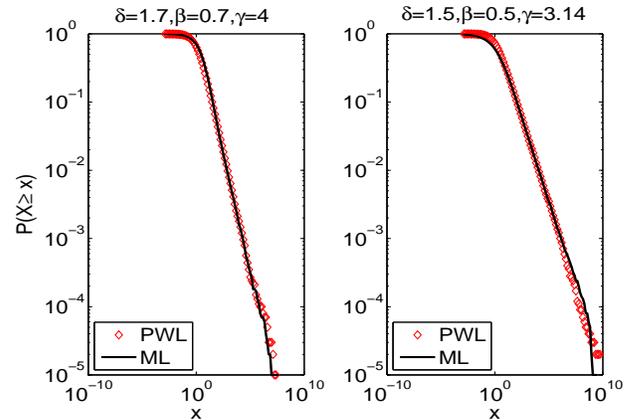}
\end{center}
\caption{(Color online) The figure shows the c.c.d.f.~of a Pareto distribution with tail exponent $\delta$ and the matching with the corresponding Mittag-Leffler$_{\gamma}(\beta)$. \textcolor{black}{The left panel presents the case $\beta=0.7$, $\gamma=4$ and $\delta=1.7$, whereas the right panel presents the case $\beta=0.5$, $\gamma=3.14$ and $\delta=1.5$. The Pareto distribution is drawn using diamonds whereas the Mittag-Leffler is drawn using a solid line.} }
\label{fig:match}

\end{figure}

On a double logarithmic scale the c.c.d.fs have a linear behavior at infinity with slope $1 -\delta= -\beta$.
Our simulations have a time horizon of $T=2000$ so the scaling $\gamma$ is chosen so that the c.c.d.f's agree well for values around and before the time horizon. The initial value of $\gamma$ to be tested for matching is the solution to the equation 
\[
	 \frac{\sin(\beta \pi)}{\pi} \frac{\Gamma(\beta)}{(t/\gamma)^\beta} = \frac{1}{t^{\delta-1}} \Longleftrightarrow \gamma  = \left(\frac{\pi}{\sin(\beta \pi)\Gamma(\beta)}\right)^{1/\beta},
\]
that implies the agreement of the asymptotic behavior of the survival functions. This first $\gamma$ choice will need to be adjusted, depending on our choice of time horizon, but the match can be achieved relatively well for moderate $\beta$ values ($\beta < 0.9$) (see Figure \ref{fig:match}).  

%\subsection{Asymptotic behaviour of the Mittag - Leffler distribution.}
%
%For the asymptotic behaviour of the c.c.d.f.~near $0$ for $|t^{\beta}| \ll 1$, one has
%\[
%E_\beta (-t^\beta) \simeq 1 - \frac{t^\beta}{\Gamma(1+\beta)}. \eqno{(a)}
%\]
%This approximation works for $t \ll 1$ with rigorous estimates coming for Taylor's theorem. Naturally this is 
%not per say a stretched exponential function, however the the expansion of a stretch exponential does have 
%the same first order Taylor expansion, 
%$\exp(-t^\beta/\Gamma(1+\beta)) \simeq1 - \frac{t^\beta}{\Gamma(1+\beta)}$.  
%We only emphasise this in order to further demonstrate remnants of the 
%Markovian behaviour that are valid for $\beta =1$. This gives the first part in equation \eqref{eq:ref:stupid}.  
%
%As an example, let us consider the case $\beta = 0.7$. 
%Fig. \ref{fig:ML-exp} shows 
%a comparison between the Mittag-Leffler function, 
%the stretched exponential defined above and the power law of equation \eqref{eq:ref:stupid}.
%\begin{figure}\label{fig:ML-exp}
%\begin{center}
%\includegraphics[height=5cm,width=7.5cm]{example.pdf}
%\end{center}
%\caption{(Color online) The Mittag-Leffler function $E_\beta (-t^\beta)$ is represented by the blue line. 
%The stretched exponential given by $\exp(-t^\beta/\Gamma(1+\beta))$ is given by the red line. 
%Finally, the power-law tail $\sin(\beta \pi) \Gamma(\beta) / (\pi t^\beta)$ is the green line.}
%\end{figure}


\bibliography{EIN_ES.bib}

%\bibliography{NMNPRL.bbl}

\end{document}